\documentclass[letterpaper, 10 pt, conference]{ieeeconf}  

\IEEEoverridecommandlockouts                              
\usepackage{graphicx} 
\usepackage{times} 

\usepackage{amsfonts,amsmath,amssymb} 

\usepackage{array}
\usepackage{subfigure}

\graphicspath{{Figures/}}

\newcommand{\HH}{{\mathcal H}}
\newcommand{\II}{{\mathcal I}}
\newcommand{\KK}{{\mathcal K}}
\newcommand{\Sal}{{\mathcal S}}

\newcommand{\OO}{{\mathcal O}}

\newcommand{\w}{{\omega}}
\newcommand{\W}{{\Omega}}
\newcommand{\eps}{{\varepsilon}}

\newcommand{\abs}[1]{\left\lvert#1\right\rvert}
\newcommand{\norm}[1]{\left\lVert#1\right\rVert}

\title{\LARGE \bf Signal injection and averaging for position estimation of Permanent-Magnet Synchronous Motors}

\author{Al~Kassem~Jebai, Fran\c{c}ois~Malrait, Philippe~Martin and~Pierre~Rouchon
\thanks{A-K.~Jebai, P.~Martin and P.~Rouchon are with the Centre Automatique et Syst\`{e}mes, MINES ParisTech, 75006~Paris,~France
{\tt\footnotesize \{al-kassem.jebai, philippe.martin, pierre.rouchon\}@mines-paristech.fr}}%
\thanks{{F.~Malrait is with Schneider Toshiba~Inverter~Europe}, 27120 Pacy-sur-Eure, France
{\tt\footnotesize francois.malrait@schneider-electric.com}}%
}

\begin{document}

\maketitle
\thispagestyle{empty}
\pagestyle{empty}

\begin{abstract}
Sensorless control of Permanent-Magnet Synchronous Motors at low velocity remains a challenging task. A now well-established method consists in injecting a high-frequency signal and use the rotor saliency, both geometric and magnetic-saturation induced. This paper proposes a clear and original analysis based on second-order averaging of how to recover the position information from signal injection; this analysis blends well with a general model of magnetic saturation. It also experimentally demonstrates the relevance for position estimation of a simple parametric saturation model recently introduced by the authors.
\end{abstract}


\section{Introduction}
Permanent-Magnet Synchronous Motors (PMSM) are widely used in industry. In the so-called ``sensorless'' mode of operation, the rotor position and velocity are not measured and the control law must make do with only current measurements. While sensorless control at medium to high velocities is well understood, with many reported control schemes and industrial products, sensorless control at low velocity remains a challenging task. The reason is that observability degenerates at zero velocity, causing a serious problem in the necessary rotor position estimation.

A now well-established method to overcome this problem is to add some persistent excitation by injecting a high-frequency signal~\cite{JanseL1995ITIA} and use the rotor saliency, whether geometric for Interior Permanent-Magnet machines~\cite{OgasaA1998ITIA,CorleL1998ITIA,AiharTYME1999ITPE}, or induced by main flux saturation for Surface Permanent-Magnet machines~\cite{JangSHIS2003ITIA}. To get a good position estimation under high-load condition it is furthermore necessary to take into account cross-saturation~\cite{GugliPV2006ITIA,BiancB2005IAC,LiZHBS2009ITIA}.

The contribution of this paper is twofold: on the one hand it proposes a clear and original analysis based on second-order averaging of how to recover the position information from signal injection; this analysis blends well with a general model of magnetic saturation including cross-saturation. On the other hand the paper experimentally demonstrates the relevance for position estimation of the saturation model introduced by the authors~\cite{JebaiMMR2011IEMDC}.

The paper is organized as follows: section~\ref{sec:model} extends the saturation model introduced in~\cite{JebaiMMR2011IEMDC} and studies its first-order observability; in section~\ref{sec:pos} position estimation by signal injection is studied thanks to second-order averaging; finally section~\ref{sec:experiment} experimentally demonstrates the relevance of the approach and the necessity of considering saturation to correctly estimate the position. 

\section{An energy-based model for the saturated PMSM}\label{sec:model}

\subsection{Notations}
In the sequel we denote by $x_{ij}:=(x_i,x_j)^T$ the vector made from the real numbers $x_i$ and $x_j$, where
$ij$ can be $dq$, $\alpha\beta$ or~$\gamma\delta$. We also define the matrices
\[M_\mu:=\begin{pmatrix}\cos\mu& -\sin\mu\\ \sin\mu& \cos\mu\end{pmatrix}
\quad\text{and}\quad
\KK:=\begin{pmatrix}0&-1\\ 1&0\end{pmatrix}, \]
and we have the useful relation
\[\frac{dM_\mu}{d\mu}=\KK M_\mu=M_\mu\KK. \]

\subsection{Energy-based model}
The model of a two-axis PMSM expressed in the synchronous $d-q$ frame reads
\begin{align}
\label{eq:dqsys1}\frac{d\phi_{dq}}{dt} &=u_{dq}-Ri_{dq}-\w\KK(\phi_{dq}+\phi_{m})\\
\frac{J}{n^2}\frac{d\omega}{dt} &= \frac{3}{2}i_{dq}^T\KK(\phi_{dq}+\phi_{m}) - \frac{\tau_L}{n}\\
\label{eq:dqsys4}\frac{d\theta}{dt} &=\omega,
\end{align}
with $\phi_{dq}$ flux linkage due to the current; $\phi_{m}:=(\lambda,0)^T$ constant flux linkage due to the permanent magnet; $u_{dq}$ impressed voltage and $i_{dq}$ stator current; $\omega$ and $\theta$ rotor (electrical) speed and position; $R$ stator resistance; $n$ number of pole pairs; $J$ inertia moment and $\tau_L$ load torque. The physically impressed voltages are $u_{\alpha\beta}:=M_\theta u_{dq}$ while the physically measurable currents are~$i_{\alpha\beta}:=M_\theta i_{dq}$.
The current can be expressed in function of the flux linkage thanks to a suitable energy function~$\HH(\phi_d,\phi_q)$ by
\begin{align}
   \label{eq:CurrentFlux}i_{dq} =\II_{dq}(\phi_{dq})
   :=\begin{pmatrix}\partial_1\HH(\phi_d,\phi_q)\\\partial_2\HH(\phi_d,\phi_q)\end{pmatrix},
\end{align}
where $\partial_k\HH$ denotes the partial derivative w.r.t. the $k^{th}$ variable~\cite{BasicMR2010ITAC,BasicJMMR2011LNCIS}; without loss of generality $\HH(0,0)=0$. Such a relation between flux linkage and current naturally encompasses cross-saturation effects.

For an unsaturated PMSM this energy function reads
\begin{align*}
\HH_l(\phi_d,\phi_q) &=\frac{1}{2L_d}\phi^2_d+\frac{1}{2L_q}\phi^2_q
\end{align*}
where $L_d$ and $L_q$ are the motor self-inductances, and we recover the usual linear relations
\begin{align*}
    i_d &=\partial_1\HH_l(\phi_d,\phi_q) =\frac{\phi_d}{L_d}\\
    i_q &=\partial_2\HH_l(\phi_d,\phi_q) =\frac{\phi_q}{L_q}.
\end{align*}

Notice the expression for $\HH$ should respect the symmetry of the PMSM w.r.t the direct axis, i.e.
\begin{equation}\label{eq:sym}
    \HH(\phi_d,-\phi_q)=\HH(\phi_d,\phi_q),
\end{equation}
which is obviously the case for~$\HH_l$.
Indeed \eqref{eq:dqsys1}--\eqref{eq:dqsys4} is left unchanged by the transformation
\begin{multline*}
(u_d,u_q,\phi_d,\phi_q,i_d,i_q,\w,\theta,\tau_L)\rightarrow\\
(u_d,-u_q,\phi_d,-\phi_q,i_d,-i_q,-\w,-\theta,-\tau_L).
\end{multline*}

\subsection{Parametric description of magnetic saturation}\label{sec:paramsat}
Magnetic saturation can be accounted for by considering a more complicated magnetic energy function~$\HH$, having $\HH_l$ for quadratic part but including also higher-order terms. From experiments saturation effects are well captured by considering only third- and fourth-order terms, hence
\begin{multline*}
\HH(\phi_d,\phi_q)=\HH_l(\phi_d,\phi_q)\\
+\sum_{i=0}^3\alpha_{3-i,i}\phi_d^{3-i}\phi_q^i+\sum_{i=0}^4\alpha_{4-i,i}\phi_d^{4-i}\phi_q^i.
\end{multline*}
This is a perturbative model where the higher-order terms appear as corrections of the dominant term~$\HH_l$.
The nine coefficients $\alpha_{ij}$
together with $L_d$, $L_q$ are motor dependent. But \eqref{eq:sym} implies $\alpha_{2,1}=\alpha_{0,3}=\alpha_{3,1}=\alpha_{1,3}=0$, so that the energy function eventually reads
\begin{multline}\label{eq:EnerSat}
\HH(\phi_d,\phi_q) =\HH_l(\phi_d,\phi_q) +\alpha_{3,0}\phi_d^3+\alpha_{1,2}\phi_d\phi_q^2 \\
+\alpha_{4,0}\phi_d^4+\alpha_{2,2}\phi_d^2\phi_q^2+\alpha_{0,4}\phi_q^4.
\end{multline}
From~\eqref{eq:CurrentFlux} and~\eqref{eq:EnerSat} the currents are then explicitly given by
\begin{align}
 \label{eq:id} i_d &=\frac{\phi_d}{L_d}+3\alpha_{3,0}\phi_d^2+\alpha_{1,2}\phi_q^2 +4\alpha_{4,0}\phi_d^3+2\alpha_{2,2}\phi_d\phi_q^2\\
 \label{eq:iq} i_q &=\frac{\phi_q}{L_q}+2\alpha_{1,2}\phi_d\phi_q+2\alpha_{2,2}\phi_d^2\phi_q+4\alpha_{0,4}\phi_q^3,
\end{align}
which are the so-called flux-current magnetization curves.


To conclude, the model of the saturated PMSM is given by~\eqref{eq:dqsys1}--\eqref{eq:dqsys4} and~\eqref{eq:id}-\eqref{eq:iq}, with $\phi_d,\phi_q,\w,\theta$ as state variables. The magnetic saturation effects are represented by the five parameters $\alpha_{3,0},\alpha_{1,2},\alpha_{4,0},\alpha_{2,2},\alpha_{0,4}$.

\subsection{Model with $i_d,i_q$ as state variables}
The model of the PMSM is usually expressed with currents as state variables. This can be achieved here by time differentiating $i_{dq}=\II_{dq}(\phi_{dq})$,
\begin{align*}
    \frac{di_{dq}}{dt}
    =D\II_{dq}(\phi_{dq})\frac{d\phi_{dq}}{dt},
\end{align*}
with $\frac{d\phi_{dq}}{dt}$ given by~\eqref{eq:dqsys1}. Fluxes are then expressed as $\phi_{dq}=\II_{dq}^{-1}(i_{dq})$ by inverting the nonlinear relations~\eqref{eq:id}-\eqref{eq:iq}; rather than performing the exact inversion, we can take advantage of the fact the coefficients $\alpha_{i,j}$ are experimentally small. At first order w.r.t. the~$\alpha_{i,j}$ we have $\phi_d=L_di_d+\OO(\abs{\alpha_{i,j}})$ and $\phi_q=L_qi_q+\OO(\abs{\alpha_{i,j}})$; plugging these expressions into~\eqref{eq:id}-\eqref{eq:iq} and neglecting $\OO(\abs{\alpha_{i,j}}^2)$ terms, we easily find
\begin{align}
 \notag\phi_d &=L_d\bigl(i_d-3\alpha_{3,0}L_d^2i_d^2 -\alpha_{1,2}L_q^2i_q^2\\
 \label{eq:phid}&\quad\qquad-4\alpha_{4,0}L_d^3i_d^3-2\alpha_{2,2}L_dL_q^2i_di_q^2\bigr)\\
 \notag\phi_q &=L_q\bigl(i_q-2\alpha_{1,2}L_dL_qi_di_q-\\
 \label{eq:phiq}&\quad\qquad2\alpha_{2,2}L_d^2L_qi_d^2i_q-4\alpha_{0,4}L_q^3i_q^3\bigr).
\end{align}

Notice the matrix
\begin{equation}\label{eq:Gmat}
\begin{pmatrix}G_{dd}(i_{dq})&G_{dq}(i_{dq})\\G_{dq}(i_{dq})&G_{qq}(i_{dq}\end{pmatrix}
:=D\II_{dq}\bigl(\II_{dq}^{-1}(i_{dq})\bigr),
\end{equation}
with coefficients easily found to be
\begin{align*}
    G_{dd}(i_{dq})&=\frac{1}{L_d}+6\alpha_{3,0}L_di_d+12\alpha_{4,0}L_d^2i_d^2+2\alpha_{2,2}L_q^2i_q^2\\
    G_{dq}(i_{dq})&=2\alpha_{1,2}L_qi_q+4\alpha_{2,2}L_di_dL_qi_q\\
    G_{qq}(i_{dq})&=\frac{1}{L_q}+2\alpha_{1,2}L_di_d+2\alpha_{2,2}L_d^2i_d^2+12\alpha_{0,4}L_q^2i_q^2,
\end{align*}
is by construction symmetric; indeed
\begin{align*}
    D\mathcal{I}_{dq}(\phi_{dq})
    =\begin{pmatrix}\partial_{11}\HH(\phi_d,\phi_q) & \partial_{21}\HH(\phi_d,\phi_q)\\
    \partial_{12}\HH(\phi_d,\phi_q) & \partial_{22}\HH(\phi_d,\phi_q)\end{pmatrix}
\end{align*}
and $\partial_{12}\HH=\partial_{21}\HH$. Therefore the inductance matrix
\begin{align*}
   \begin{pmatrix}L_{dd}(i_{dq}) & L_{dq}(i_{dq})\\ L_{dq}(i_{dq})& L_{qq}(i_{dq})\end{pmatrix}
   :=&\begin{pmatrix}G_{dd}(i_{dq})&G_{dq}(i_{dq})\\G_{dq}(i_{dq})&G_{qq}(i_{dq})\end{pmatrix}^{-1}.
\end{align*}
is also symmetric, though this is not always acknowledged in saturation models encountered in the literature.

\subsection{First-order observability}
To emphasize position estimation is not easy at low speed, we study the observability of~\eqref{eq:dqsys1}--\eqref{eq:dqsys4}, augmented by
\begin{align}\label{eq:constantload}
\frac{d\tau_L}{dt}=0,
\end{align}
around a permanent trajectory defined by
\begin{align*}
0 &=\overline u_{dq}-R\overline i_{dq}-\w\KK(\overline\phi_{dq}+\phi_{m})\\
0 &=\frac{3}{2}\overline i_{dq}^T\KK(\overline\phi_{dq}+\phi_{m}) - \frac{\overline\tau_L}{n}\\
\frac{d\overline\theta}{dt} &=\overline\omega.
\end{align*}
Notice such a permanent trajectory is not a steady state point unless $\overline\w:=0$ since $\overline\theta$ hence $u_{\alpha\beta}$ and $i_{\alpha\beta}$ are time-varying ($\overline\phi_{dq},\overline i_{dq},\overline u_{dq},\overline\w,\overline\tau_L$ are on the other hand constant). Also the rationale for augmenting the model with~\eqref{eq:constantload} is that it is usually also desired to estimate an unknown load torque.

From $i_{\alpha\beta}=M_\theta i_{dq}$ we then get
\begin{align}\label{eq:deltai}
\delta i_{\alpha\beta}=\delta M_\theta\overline i_{dq}+M_{\overline\theta}\delta i_{dq}
=M_{\overline\theta}(\KK\overline i_{dq}\delta\theta+\delta i_{dq}),
\end{align}
and similarly for~$u_{\alpha\beta}$. The linearization of~\eqref{eq:dqsys1}--\eqref{eq:dqsys4} and~\eqref{eq:constantload} around a permanent trajectory is then
\begin{align*}
\frac{d\delta\phi_{dq}}{dt} &= M^T_{\overline\theta}(\delta u_{\alpha\beta}-R\delta i_{\alpha\beta})
-\KK(\overline\phi_{dq}+\phi_{m})\delta\w\\
&\quad+\overline\w\bigl((\overline\phi_{dq}+\phi_{m})\delta\theta-\KK\delta\phi_{dq}\bigr)\\
\frac{J}{n^2}\frac{d\delta\omega}{dt} &=
\frac{3}{2}\delta i_{\alpha\beta}^TM_{\overline\theta}\KK(\overline\phi_{dq}+\phi_{m})\\
&\quad-\frac{3}{2}\overline i_{dq}^T\bigl((\overline\phi_{dq}+\phi_{m})\delta\theta-\KK\delta\phi_{dq}\bigr)
-\frac{\delta\tau_L}{n}\\
\frac{d\delta\theta}{dt} &=\delta\omega\\
\frac{d\delta\tau_L}{dt} &=0,
\end{align*}
where we have used $\overline u_{dq}-R\overline i_{dq}=\overline\w\KK(\overline\phi_{dq}+\phi_{m})$.

On the other hand time differentiating $\delta i_{\alpha\beta}$ yields
\begin{align*}
\frac{d\delta i_{\alpha\beta}}{dt} &=\frac{dM_{\overline\theta}}{dt}(\KK\overline i_{dq}\delta\theta+\delta i_{dq})
+M_{\overline\theta}\Bigl(\KK\overline i_{dq}\frac{d\delta\theta}{dt}+\frac{d\delta i_{dq}}{dt}\Bigr)\\
&=\overline\w\KK\delta i_{\alpha\beta}+M_{\overline\theta}\Bigl(\KK\overline i_{dq}\delta\w
+D\II_{dq}(\overline\phi_{dq})\frac{d\delta\phi_{dq}}{dt}\Bigr)
\end{align*}
where we have used $\delta i_{dq}=D\II_{dq}(\overline\phi_{dq})\delta\phi_{dq}$. Therefore
\begin{align*}
&[D\II_{dq}(\overline\phi_{dq})]^{-1}M_{\overline\theta}^T\frac{d\delta i_{\alpha\beta}}{dt}\\
&=
[D\II_{dq}(\overline\phi_{dq})]^{-1}\KK\overline i_{dq}\delta\w+\frac{d\delta\phi_{dq}}{dt}
+\text{LC}(\delta i_{\alpha\beta},\delta u_{\alpha\beta})\\
&=\KK\Phi(\overline\phi_{dq})\delta\w
+\Phi(\overline\phi_{dq})\overline\w\delta\theta
+\text{LC}(\delta i_{\alpha\beta},\delta u_{\alpha\beta})
\end{align*}
where
$$\Phi(\overline\phi_{dq}):=\phi_{m}+\overline\phi_{dq}
+\KK[D\II_{dq}(\overline\phi_{dq})]^{-1}\KK\II_{dq}(\overline\phi_{dq})$$
and $\text{LC}(\delta i_{\alpha\beta},\delta u_{\alpha\beta})$ is some matrix linear combination of $\delta i_{\alpha\beta}$ and~$\delta u_{\alpha\beta}$.
Similarly
\begin{align}\label{eq:omegaobs}
\frac{J}{n^2}\frac{d\delta\omega}{dt}=-\frac{3}{2}\Phi(\overline\phi_{dq})\delta\theta -\frac{\delta\tau_L}{n}+\text{LC}(\delta i_{\alpha\beta},\delta u_{\alpha\beta}).
\end{align}

Now $\Phi(\overline\phi_{dq})$ is non-zero in any reasonable situation, hence $\Phi(\overline\phi_{dq})$ and $\KK\Phi(\overline\phi_{dq})$ are independent vectors. If $\overline\w\neq0$ it is thus clear that $\delta\theta$ and $\delta\w$ can be expressed in function of $i_{\alpha\beta},\frac{d\delta i_{\alpha\beta}}{dt}$ and~$u_{\alpha\beta}$; as a consequence $\delta\tau_L$ is a function of $i_{\alpha\beta},\frac{d\delta i_{\alpha\beta}}{dt},\frac{d^2\delta i_{\alpha\beta}}{dt^2}$ and~$u_{\alpha\beta},\frac{d\delta u_{\alpha\beta}}{dt}$ by~\eqref{eq:omegaobs}; finally $\delta\phi_{dq}$ is by~\eqref{eq:deltai} a function of~$\delta i_{\alpha\beta}$. In other words the linearized system is observable by~$\delta i_{\alpha\beta}$.

If $\overline\w=0$ only $\delta\w$ can be recovered from $i_{\alpha\beta},\frac{d\delta i_{\alpha\beta}}{dt}$ and~$u_{\alpha\beta}$; as a consequence only $\Phi(\overline\phi_{dq})\delta\theta+\frac{\delta\tau_L}{n}$ can be recovered from~\eqref{eq:omegaobs}, and nothing new is gained by further differentiation. In other words the linearized system is not observable, as pointed out in~\cite{BasicMR2010ITAC}.

Estimating the rotor position at low speed is thus inherently difficult. Yet it is doable when the motor exhibits some saliency, geometric or saturation-induced, with the help of permanent excitation such as signal injection. Indeed a more thorough analysis along the lines of~\cite{VaclaB2007CCA,ZaltnAGB2009SCS} would reveal the system is in that case observable in the nonlinear sense provide the motor exhibits some saliency.

\section{Position estimation by high frequency voltage injection}\label{sec:pos}

\subsection{Signal injection and averaging}
A general sensorless control law can be expressed as
\begin{align}
    \label{eq:claw1}u_{\alpha\beta} &=M_{\theta_c}u_{\gamma\delta}\\
    \frac{d\theta_c}{dt} &=\omega_c\\
    \frac{d\eta}{dt} &=a\bigl(M_{\theta_c}i_{\gamma\delta},\theta_c,\eta,t\bigr)\\
    \label{eq:claw4}\omega_c &=\W_c\bigl(M_{\theta_c}i_{\gamma\delta},\theta_c,\eta,t\bigr)\\
    \label{eq:claw5}u_{\gamma\delta} &=\mathcal U_{\gamma\delta}\bigl(M_{\theta_c}i_{\gamma\delta},\theta_c,\eta,t\bigr),
\end{align}
where the measured currents $i_{\alpha\beta}=M_{\theta_c}i_{\gamma\delta}$ are used to compute $u_{\gamma\delta}$, $\omega_c$ and the evolution of the internal (vector) variable $\eta$ of the controller; $\theta_c$ and~$\omega_c$ are known by design.

It will be convenient to write the system equations~\eqref{eq:dqsys1}--\eqref{eq:dqsys4} in the $\gamma-\delta$ frame defined by $x_{\gamma\delta}:=M_{\theta-\theta_c}x_{dq}$, which gives
\begin{align}
\label{eq:gdsys1}\frac{d\phi_{\gamma\delta}}{dt}
&=u_{\gamma\delta}-Ri_{\gamma\delta}-\w_c\KK\phi_{\gamma\delta}
-\w\KK M_{\theta-\theta_c}\phi_{m}\\
\frac{J}{n^2}\frac{d\omega}{dt} &= \frac{3}{2}i_{\gamma\delta}^T\KK(\phi_{\gamma\delta}+M_{\theta-\theta_c}\phi_{m}) - \frac{\tau_L}{n}\\
\label{eq:gdsys4}\frac{d\theta}{dt} &=\omega,
\end{align}
where from~\eqref{eq:id}-\eqref{eq:iq} currents and fluxes are related by
\begin{equation}\label{eq:iphigd}
    i_{\gamma\delta}=M_{\theta-\theta_c}\mathcal{I}_{dq}(M^T_{\theta-\theta_c}\phi_{\gamma\delta}).
\end{equation}

To estimate the position we will superimpose on some desirable control law~\eqref{eq:claw5} a fast-varying pulsating voltage,
\begin{equation}\label{eq:HFvoltage}
    u_{\gamma\delta} =\mathcal U_{\gamma\delta}\bigl(M_{\theta_c}i_{\gamma\delta},\theta_c,\eta,t\bigr)
     + \widetilde u_{\gamma\delta}f(\Omega t),
\end{equation}
where $f$ is a $2\pi$-periodic function with zero mean and $\widetilde u_{\gamma\delta}$ could like~$\cal U_{\gamma\delta}$ depend on $M_{\theta_c}i_{\gamma\delta},\theta_c,\eta,t$ (though it is always taken constant in the sequel). The constant pulsation $\Omega$ is chosen ``large'', so that $f(\Omega t)$ can be seen as a ``fast'' oscillation; typically $\Omega:=2\pi\times500~\text{rad/s}$ in the experiments in section~\ref{sec:experiment}.

If we apply this modified control law to~\eqref{eq:gdsys1}--\eqref{eq:gdsys4}, then it can be shown the solution of the closed loop system is
\begin{align}
    \label{eq:gdavrth1}\phi_{\gamma\delta} &= \overline\phi_{\gamma\delta} + \frac{\widetilde  u_{\gamma\delta}}{\Omega}F(\Omega t) + \OO{(\frac{1}{\Omega^2})}\\
    \omega &= \overline\omega + \OO{(\frac{1}{\Omega^2})} \\
    \theta &= \overline\theta + \OO{(\frac{1}{\Omega^2})} \\
    \theta_c &= \overline\theta_c + \OO{(\frac{1}{\Omega^2})} \\
    \label{eq:gdavrth5}\eta &=  \overline\eta + \OO{(\frac{1}{\Omega^2})},
\end{align}
where $F$ is the primitive of $f$ with zero mean ($F$ clearly has the same period as~$f$); $(\overline\phi_{\gamma\delta},\overline\omega,\overline\theta,\overline\theta_c,\overline\eta)$ is the ``slowly-varying'' component of $(\phi_{\gamma\delta},\omega,\theta,\theta_c,\eta)$, i.e. satisfies
\begin{align*}
\frac{d\overline\phi_{\gamma\delta}}{dt}
&=\overline u_{\gamma\delta}-R\overline i_{\gamma\delta}-\overline\w_c\KK\overline\phi_{\gamma\delta}
-\w\KK M_{\overline\theta-\overline\theta_c}\phi_{m}\\
\frac{J}{n^2}\frac{d\overline\omega}{dt} &=\frac{3}{2}\overline i_{\gamma\delta}^T\KK(\overline\phi_{\gamma\delta}+M_{\overline\theta-\overline\theta_c}\phi_{m})
-\frac{\tau_L}{n}\\
    \frac{d\overline\theta}{dt}& = \overline\omega\\
    \frac{d\overline\theta_c}{dt} &=\overline\omega_c\\
    \frac{d\overline\eta}{dt} &=a\bigl(M_{\overline\theta_c}\overline i_{\gamma\delta},\overline\theta_c,\overline\eta,t\bigr),
\end{align*}
where
\begin{align}
    \label{eq:bari}\overline i_{\gamma\delta} &=M_{\overline\theta-\overline\theta_c}\mathcal{I}_{dq}(M^T_{\overline\theta-\overline\theta_c}\overline\phi_{\gamma\delta})\\
    \notag\overline\omega_c &=\W_c\bigl(M_{\overline\theta_c}\overline i_{\gamma\delta},\overline\theta_c,\overline\eta,t\bigr)\\
    \notag\overline u_{\gamma\delta} &=\mathcal U_{\gamma\delta}\bigl(M_{\overline\theta_c}\overline i_{\gamma\delta},\overline\theta_c,\overline\eta,t\bigr).
\end{align}
Notice this slowly-varying system is exactly the same as~\eqref{eq:gdsys1}--\eqref{eq:gdsys4} acted upon by the unmodified control law~\eqref{eq:claw1}--\eqref{eq:claw5}. In other words adding signal injection:
\begin{itemize}
  \item has a very small effect of order $\OO{(\frac{1}{\Omega^2})}$ on the mechanical variables $\theta,\omega$ and the controller variables $\theta_c,\eta$
  \item has a small effect of order $\OO{(\frac{1}{\Omega})}$ on the flux~$\phi_{\gamma\delta}$; this effect will be used in the next section to extract the position information from the measured currents.
\end{itemize}

The proof relies on a direct application of second-order averaging of differential equations, see~\cite{SandeVM2007book} section~$2.9.1$ and for the slow-time dependance section~$3.3$. Indeed setting $\eps:=\frac{1}{\W}$ and $x:=(\phi_{\gamma\delta},\omega,\theta,\theta_c,\eta)$, \eqref{eq:gdsys1}--\eqref{eq:gdsys4} acted upon by the modified control law \eqref{eq:claw1}--\eqref{eq:claw4} and~\eqref{eq:HFvoltage} is in the so-called standard form for averaging (with slow-time dependance)
\begin{align*}
    \frac{dx}{d\sigma} &=\eps f_1(x,\eps\sigma,\sigma)
    :=\eps\bigl(\overline f_1(x,\eps\sigma)+\widetilde f_1(x,\eps\sigma)f(\sigma)\bigr),
\end{align*}
with $f_1$ $T$-periodic w.r.t. its third variable and $\eps$ as a small parameter.  Therefore its solution can be approximated as 
\begin{align*}
    x(\sigma)&=z(\sigma)+\eps\bigl(u_1(z(\sigma,\eps\sigma,\sigma)\bigr)+\OO(\eps^2),
\end{align*}
where $z(\sigma)$ is the solution of
\begin{align*}
    \frac{dz}{d\sigma} &=\eps g_1(z,\eps\sigma) + \eps^2g_2(z,\eps\sigma)
\end{align*}
and
\begin{align*}
    g_1(y,\eps\sigma) &:=\frac{1}{T}\int_0^Tf_1(y,\eps\sigma,s)ds
    =\overline f_1(y,\eps\sigma)\\
    v_1(y,\eps\sigma,\sigma) &:=\int_0^\sigma\bigl(f_1(y,\eps\sigma,s)-g_1(y,\eps\sigma)\bigr)ds\\
    &=\widetilde f_1(y,\eps\sigma)\int_0^\sigma f(s)ds\\
    u_1(y,\eps\sigma,\sigma) &:=v_1(y,\eps\sigma,\sigma)-\frac{1}{T}\int_0^Tv_1(y,\eps\sigma,s)ds\\
    &=\widetilde f_1(y,\eps\sigma)F(\sigma)\\
    K_2(y,\eps\sigma,\sigma)\\:= \partial_1f_1(y,&\eps\sigma,\sigma)u_1(y,\eps\sigma,\sigma)
    -\partial_1u_1(y,\eps\sigma,\sigma)g_1(y,\eps\sigma)\\
    =[\overline f_1,\widetilde f_1](&y,\eps\sigma)F(\sigma)
    +\frac{1}{2}\partial_1\widetilde f_1(y,\eps\sigma)\widetilde f_1(y,\eps\sigma)
    \frac{dF^2(\sigma)}{d\sigma}\displaybreak[3]\\
    g_2(y,\eps\sigma) &:=\frac{1}{T}\int_0^TK_2(y,\eps\sigma,s)ds=0.
\end{align*}
We have set
\[[\overline f_1,\widetilde f_1](y,\eps\sigma)\!:=\!\partial_1\overline f_1(y,\eps\sigma)\widetilde f_1(y,\eps\sigma)-\partial_1\widetilde f_1(y,\eps\sigma)\overline f_1(y,\eps\sigma)\]
and $F(\sigma):=\int_0^\sigma f(s)ds-\frac{1}{T}\int_0^T\int_0^\sigma f(s)dsd\sigma$, i.e. $F$ is the (of course $T$-periodic) primitive of~$f$ with zero mean.

Translating back to the original variables eventually yields the desired result~\eqref{eq:gdavrth1}--\eqref{eq:gdavrth5}.

\subsection{Position estimation}\label{sec:posesti}
We now express the effect of signal injection on the currents: plugging~\eqref{eq:gdavrth1} into~\eqref{eq:iphigd} we have
\begin{align}
    \notag i_{\gamma\delta}
    &=M_{\overline\theta-\overline\theta_c+\OO{(\frac{1}{\Omega^2})}}\\
    \notag&\qquad\mathcal{I}_{dq}\Bigl(M^T_{\overline\theta-\overline\theta_c+\OO{(\frac{1}{\Omega^2})}}
    \bigl(\overline\phi_{\gamma\delta}
    + \frac{\widetilde  u_{\gamma\delta}}{\Omega}F(\Omega t) + \OO{(\frac{1}{\Omega^2})}\bigr)\Bigr)\\
    \label{eq:ilowhigh}&=\overline i_{\gamma\delta}+\widetilde i_{\gamma\delta}F(\Omega t)+\OO{(\frac{1}{\Omega^2})},
\end{align}
where we have used~\eqref{eq:bari} and performed a first-order expansion to get
\begin{align}
    \notag\widetilde i_{\gamma\delta}
    &:=M_{\overline\theta-\overline\theta_c}D\mathcal{I}_{dq}\bigl(M^T_{\overline\theta-\overline\theta_c}
    \overline\phi_{\gamma\delta}\bigr)M^T_{\overline\theta-\overline\theta_c}\frac{\widetilde  u_{\gamma\delta}}{\Omega}\\
    \label{eq:tildei}&=M_{\overline\theta-\overline\theta_c}D\mathcal{I}_{dq}
    \Bigl(\mathcal{I}_{dq}^{-1}\bigl(M^T_{\overline\theta-\overline\theta_c}
    \overline i_{\gamma\delta}\bigr)\Bigr)M^T_{\overline\theta-\overline\theta_c}\frac{\widetilde  u_{\gamma\delta}}{\Omega}.
\end{align}
We will see in the next section how to recover $\widetilde i_{\gamma\delta}$ and $\overline i_{\gamma\delta}$ from the measured currents~$i_{\gamma\delta}$.
Therefore \eqref{eq:tildei} gives two (redundant) relations relating the unknown angle~$\overline\theta$ to the known variables $\overline\theta_c,\widetilde i_{dq},\overline i_{\gamma\delta},\widetilde u_{dq}$, provided the matrix
\[ \Sal(\mu,\overline i_{\gamma\delta}):=M_\mu D\mathcal{I}_{dq}
    \Bigl(\mathcal{I}_{dq}^{-1}\bigl(M^T_\mu\overline i_{\gamma\delta}\bigr)\Bigr)M^T_\mu \]
effectively depends on its first argument~$\mu$. This ``saliency condition'' is what is needed to ensure nonlinear observability.
The explicit expression for $\Sal(\mu,\overline i_{\gamma\delta})$ is obtained thanks to~\eqref{eq:Gmat}. In the case of an unsaturated magnetic circuit this matrix boils down to
\begin{align*}
    \Sal(\mu,\overline i_{\gamma\delta}) &=M_\mu
    \begin{pmatrix}\frac{1}{L_d}&0\\0&\frac{1}{L_q}\end{pmatrix}M^T_\mu\\
    &=\textstyle\frac{L_d+L_q}{2L_dL_q}
    \begin{pmatrix}1+\frac{L_d-L_q}{L_d+L_q}\cos2\mu &\frac{L_d-L_q}{L_d+L_q}\sin2\mu\\
    \frac{L_d-L_q}{L_d+L_q}\sin2\mu &1-\frac{L_d-L_q}{L_d+L_q}\cos2\mu\end{pmatrix}
\end{align*}
and does not depend on~$i_{\gamma\delta}$; notice this matrix does not depend on~$\mu$ for an unsaturated machine with no geometric saliency. Notice also \eqref{eq:tildei} defines in that case two solutions on~$]-\pi,\pi]$ for the angle~$\overline\theta$ since $\Sal(\mu,\overline i_{\gamma\delta})$ is actually a function of~$2\mu$; in the saturated case there is generically only one solution, except for some particular values of~$\overline i_{\gamma\delta}$.

There are several ways to extract the rotor angle information from~\eqref{eq:tildei}, especially for real-time use inside a feedback law. In this paper we just want to demonstrate the validity of~\eqref{eq:tildei} and we will be content with directly solving it through a nonlinear least square problem; in other words we estimate the rotor position as
\begin{equation}\label{eq:nonlin}
      \widehat\theta = \theta_c + \arg\min_{\mu\in]-\pi,\pi]}\norm{
      \widetilde i_{\gamma\delta} - \Sal(\mu,\overline i_{\gamma\delta})\frac{\widetilde u_{\gamma\delta}}{\Omega}}^2.
\end{equation}

\subsection{Current demodulation}\label{sec:currdem}
To estimate the position information using e.g.~\eqref{eq:nonlin} it is necessary to extract the low- and high-frequency components $\overline i_{\gamma\delta}$ and $\widetilde i_{\gamma\delta}$ from the measured current~$i_{\gamma\delta}$.
Since by~\eqref{eq:ilowhigh} $i_{\gamma\delta}(t)\approx\overline i_{\gamma\delta}(t)+\widetilde i_{\gamma\delta}(t)F(\Omega t)$ with $\overline i_{\gamma\delta}$ and $\widetilde i_{\gamma\delta}$ by construction nearly constant on one period of~$F$, we may write
\begin{align*}
    \overline i_{\gamma\delta}(t) &\approx\frac{1}{T}\int_{t-T}^ti_{\gamma\delta}(s)ds\\
    \widetilde i_{\gamma\delta}(t) &
    \approx\frac{\int_{t-T}^ti_{\gamma\delta}(s)F(\Omega s)ds}{\int_0^TF^2(\Omega s)ds},
\end{align*}
where $T:=\frac{2\pi}{\Omega}$. Indeed as $F$ is $2\pi$-periodic with zero mean,
\begin{align*}
    \int_{t-T}^ti_{\gamma\delta}(s)ds &\approx
    \overline i_{\gamma\delta}(t)\int_{t-T}^tds + \widetilde i_{\gamma\delta}(t)\int_{t-T}^tF(\Omega s)ds\\
    &=T\overline i_{\gamma\delta}(t)\\
    \int_{t-T}^ti_{\gamma\delta}(s)F(\Omega s)ds &\approx
    \overline i_{\gamma\delta}(t)\int_{t-T}^tF(\Omega s)ds \\
    &\qquad+ \widetilde i_{\gamma\delta}(t)\int_{t-T}^tF^2(\Omega s)ds\\
    &=\widetilde i_{\gamma\delta}(t)\int_0^TF^2(\Omega s)ds.
\end{align*} 
\section{Experimental Results}\label{sec:experiment}
\begin{figure}[t]
\centerline{\includegraphics[width=\columnwidth]{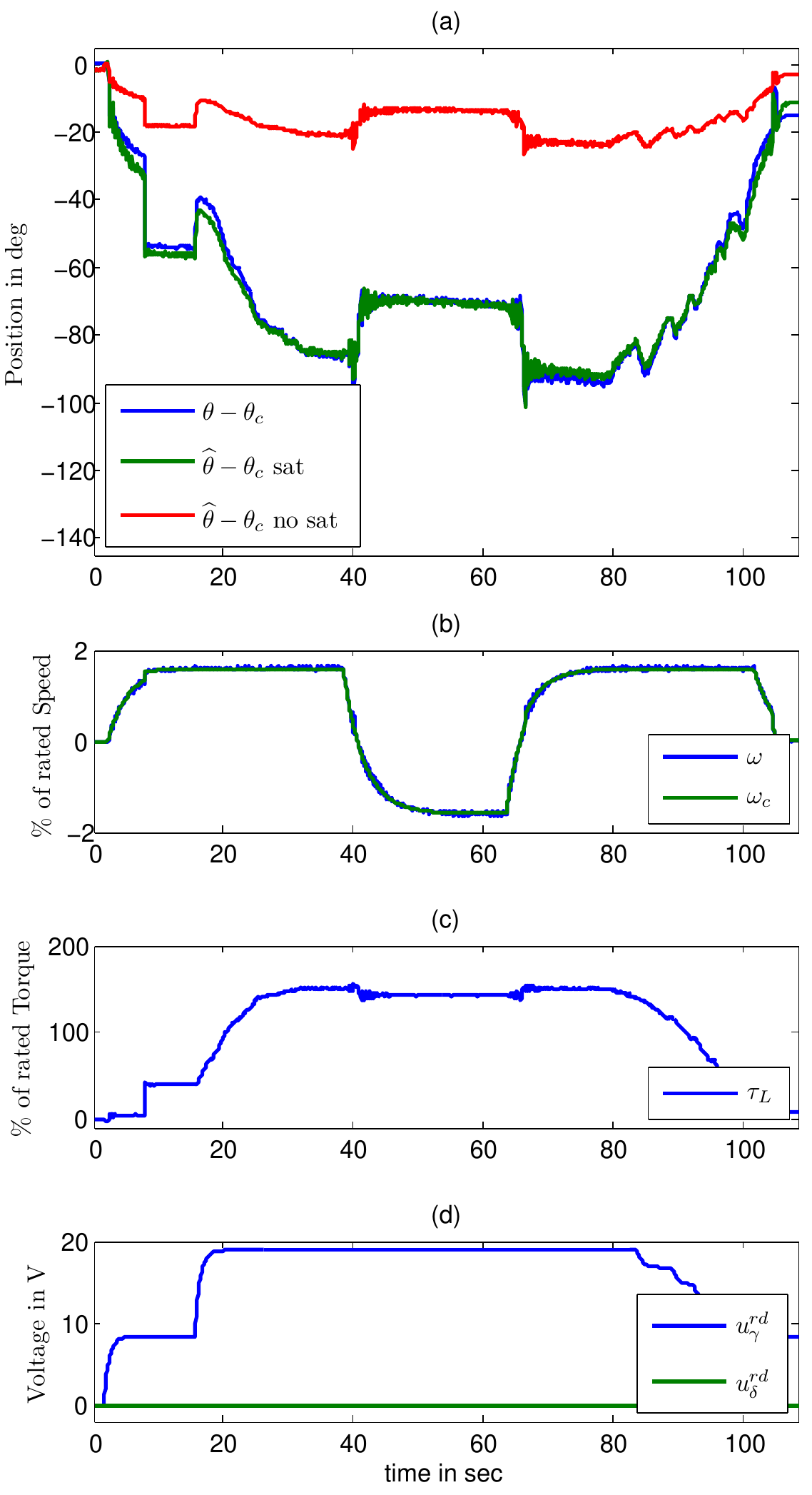}}
 \caption{Long test under various conditions:
 (a) measured $\theta-\theta_c$, estimated $\widehat\theta-\theta_c$ with and without saturation model; (b) measured speed $\omega$, reference speed~$\omega_c$; (c) load torque $\tau_L$; (d) voltages $u_{\gamma\delta}^{rd}$.}\label{fig:LongTest}
\end{figure}
\begin{table}\setlength{\extrarowheight}{2pt}
\caption{Rated and saturation parameters of the test motor.\label{tab:param}}
\begin{center}
\begin{tabular}{| c | c || c | c | }
  \firsthline
  Rated power &$1500$ W        & $L_d$& $7.9$ mH \\
  \hline
  Rated current $I_n$ (peak) & $5.19$ A     & $L_q$& $8.2$ mH\\
  \hline
  Rated speed &  $3000$~rpm       & $\alpha_{3,0} L_d^2I_n$ & $0.0551$\\
  \hline
  Rated torque & $6.06$ Nm     & $\alpha_{1,2}L_dL_qI_n$ & $0.0545$\\
  \hline
  n & $5$     & $\alpha_{4,0}L_d^3I_n^2$ & $0.0170$\\
  \hline
  R & $2.1$ $\Omega$     & $\alpha_{2,2}L_dL_q^2I_n^2$ & $0.0249$\\
   \hline
  $\lambda$  & $0.155$ $W_b$     & $\alpha_{0,4}L_q^3I_n^2$ & $0.0067$\\
   \lasthline
\end{tabular}
\end{center}
\end{table} 
\subsection{Experimental setup}
The methodology developed in the paper was tested on a surface-mounted PMSM with parameters listed in table~\ref{tab:param}. The magnetic parameters were estimated using the procedure of~\cite[section~III]{JebaiMMR2011IEMDC}. Notice this motor has little geometric saliency ($L_d\approx L_q$) hence the saturation-induced saliency is paramount to estimate the rotor position.

The experimental setup consists of an industrial inverter ($400$~V DC bus, $4$~kHz PWM frequency), an incremental encoder, a dSpace fast prototyping system with 3 boards (DS1005, DS5202 and EV1048). The measurements are sampled also at~$4$~kHz, and synchronized with the PWM frequency. The load torque is created by a $4$~kW DC motor.

\subsection{Validation of the rotor position estimation procedure}
\subsubsection{Control law}
Since the goal is only to test the validity of the angle estimation procedure, a very simple $V/f$ open-loop (i.e. $\Omega_c$ and $\mathcal U_{\gamma\delta}$ do not depend on $i_{\gamma\delta}$) control law is used for~\eqref{eq:claw1}--\eqref{eq:claw5}; a fast-varying ($\Omega:=2\pi\times500~\text{rad/s}$) square voltage with constant amplitude is added in accordance with~\eqref{eq:HFvoltage}, resulting in
\begin{align*}
    \frac{d\theta_c}{dt} &= \w_c(t)\\
       u_{\gamma\delta} &= u_{\gamma\delta}^{rd}(t) + \w_c(t)\phi_{m} + \widetilde u_{\gamma\delta} f(\Omega t).
\end{align*}
Here $\w_c(t)$ is the motor speed reference; $u_{\gamma\delta}^{rd}(t)$ is a filtered piece-wise constant vector compensating the resistive voltage drop in order to maintain the torque level and the motor stability;
$\widetilde u_{\gamma\delta}:=(\widetilde u,0)^T$ with $\widetilde u:=15$~V.

The rotor position $\widehat\theta$ is then estimated according to~\eqref{eq:nonlin}.

\subsubsection{Long test under various conditions~(Fig.~\ref{fig:LongTest})}
Speed and torque are slowly changed over a period of about two minutes; the speed remains between~$\pm2\%$ of the rated speed and the torque varies from $0\%$ to $150\%$ of the rated toque. This represents typical operation conditions at low speed.

When the saturation model is used the agreement between the estimated position $\widehat{\theta}$ and the measured position $\theta$ is very good, with an error always smaller than a few (electrical) degrees. By contrast the estimated error without using the saturation model (i.e. with all the magnetic saturation parameters $\alpha_{ij}$ taken to zero) can reach up to $70^\circ$ (electrical). This demonstrates the importance of considering an adequate saturation model.

\begin{figure}[ht]
\centerline{\includegraphics[width=\columnwidth]{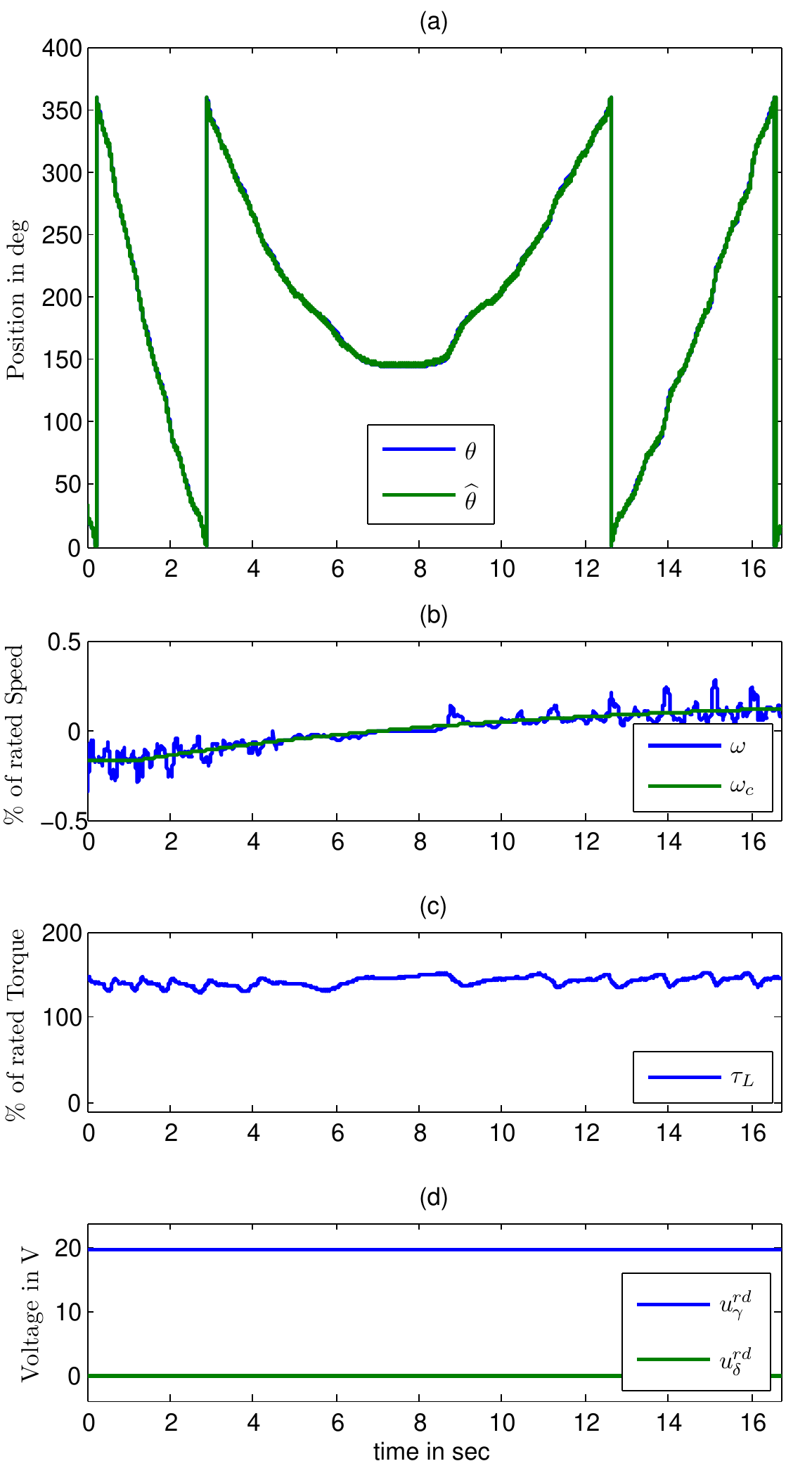}}
\caption{Very slow speed reversal:
(a) measured $\theta$, estimated $\widehat\theta$ with 
saturation model; (b) measured speed $\omega$, reference speed $\omega_c$; (c) load torque~$\tau_L$; (d) voltages $u_{\gamma\delta}^{rd}$.}\label{fig:SpeedReversal}
\end{figure}

\subsubsection{Very slow speed reversal~(Fig.~\ref{fig:SpeedReversal})}
 The speed is very slowly (in about $20$ seconds) changed from $-0.2\%$ to $+0.2\%$ of the rated speed at~$150\%$ of the rated torque. This is a very demanding test since the motor always remains in the first-order unobservability region, moreover under high load. Once again the estimated angle closely agrees with the measured angle.

\section{Conclusion}
We have presented a new procedure based on signal injection for estimating the rotor angle of a PMSM at low speed, with an original analysis based on second-order averaging. This is not an easy problem in view of the observability degenaracy at zero speed.  The method is general in the sense it can accommodate virtually any control law, saturation model, and form of injected signal. The relevance of the method and the importance of using an adequate magnetic saturation model has been experimentally demonstrated on a surface-mounted PMSM with little geometric saliency.


\bibliographystyle{IEEEtran}
\bibliography{biblio2012CDC}

%
%
%
%

\end{document}